\documentclass{article}
\usepackage{amsfonts,latexsym}
\usepackage{amsmath,amssymb}
\usepackage{amsthm}

\DeclareMathOperator{\SLtwo}{SL_2}

\newcommand{\C}{\mathbb{C}}

\newcommand{\oo}{\mathcal{O}}
\newcommand{\nn}{\mathcal{N}}
\newcommand{\vv}{\mathcal{V}}
\newcommand{\divides}{\mid}
\newcommand{\supr}[1]{^{\,#1}}

\newtheorem{theorem}{Theorem}[section]
\newtheorem{proposition}[theorem]{Proposition}
\newtheorem{lemma}[theorem]{Lemma}

\theoremstyle{definition}

\title{The invariants of the binary decimic}
\author{Andries E. Brouwer \& Mihaela Popoviciu\thanks{The second author
is partially supported by the Swiss National Science Foundation.}}
\date{July 31, 2009}

\begin{document}
\maketitle

\begin{abstract}
We consider the algebra of invariants of binary forms of degree $10$
with complex coefficients, construct a system of parameters with degrees
2, 4, 6, 6, 8, 9, 10, 14 and find the 106 basic invariants.
\end{abstract}

\section{Introduction}
{\bf Invariants}\\[1pt]
Let $\oo (V_n)^{\SLtwo}$ denote the algebra of invariants of binary
forms (forms in two variables) of degree $n$ with complex coefficients.
This algebra was extensively studied in the nineteenth century, 
and for $n \leq 6$ the structure was clear and a finite basis was known.
Gordan \cite{Go} proved in 1868 that $\oo (V_n)^{\SLtwo}$ has
a finite basis for all $n$. For $n=7$ the invariants were determined
by von Gall \cite{vG7} and Dixmier \& Lazard \cite{DiLa}
(see also Bedratyuk \cite{Be}).
The invariants for $n=8$ were found by von Gall \cite{vG8} and
Shioda \cite{Sh}.
The case $n=9$ was done by Cr\"oni \cite{Cr} and the present authors
\cite{BP}.
Here we consider the case $n=10$, and show that $\oo (V_{10})^{\SLtwo}$
is generated by 106 (explicitly known) basic invariants, and give
the degrees.

\begin{proposition}\label{prop1}
The algebra $I$ of invariants of the binary decimic (form of degree $10$)
is generated by $106$ invariants. The nonzero numbers $d_m$ of
basic invariants of degree $m$ are \\
\smallskip
\begin{tabular}{@{\,}c|cccccccccccccccc@{\,}}
$m$ & $2$ &$4$ & $6$ & $8$ & $9$ & $10$ & $11$ & $12$ & $13$ & $14$ &
$15$ & $16$ & $17$ & $18$ & $19$ & $21$ \\
\hline
$d_m$ & $1$ & $1$ & $4$ & $5$ & $5$ & $8$ & $8$ & $12$ & $15$ & $13$ &
$19$ & $5$ & $5$ & $1$ & $2$ & $2$ \\
\end{tabular}
\end{proposition}

This list agrees with Sylvester \& Franklin \cite{SyFr} for degrees
less than 17. Sylvester predicted 3 basic invariants of degree 17
and none of degree higher than 17 for a total of 99 basic invariants.
Tom Hagedorn (unpublished) found 104 invariants,
cf. Olver \cite{Olver} (p. 40).
The existence of basic invariants of degree 21 seems to be new.
That the list is complete follows as a corollary from
the construction of a homogeneous system of parameters (hsop), see below.

\medskip\noindent
{\bf Systems of parameters}\\
A (homogeneous) {\it system of parameters} for a graded algebra $A$
is an algebraically independent set $S$ of homogeneous elements of $A$
such that $A$ is module-finite over the subalgebra generated by the set $S$.
Hilbert \cite{Hi2} showed the existence of a system of parameters
for algebras of invariants, cf.~Proposition \ref{hilbert} below.

Here we find an explicit system of parameters for $\oo(V_{10})^{\SLtwo}$.

\begin{proposition}\label{hsop}
The algebra $I$ of invariants of the binary decimic has a
system of parameters of degrees $2$, $4$, $6$, $6$, $8$, $9$, $10$, $14$.
\end{proposition}

This is useful, since it provides an upper bound for the degrees
of basic invariants that is sufficiently low, so that a simple
computer search can find a basis for the invariants up to that degree.

\section{Finding the basic invariants}\label{finding}
A set of basic invariants of the algebra $I$ of invariants
is a minimal set of generators. The individual generators are
not uniquely determined, but their degrees are.

\medskip
The ring $I$ is graded: $I = \oplus_m I_m$, where $I_m$ is the
subspace of invariants, homogeneous of degree $m$.
If $y_1,\ldots,y_{n-2}$ is a system of parameters, where $y_i$ has
degree $d_i$, then the Poincar\'e series $P(t)$, defined by
$P(t) = \sum_m \dim I_m t^m$, can be written as a rational function
in $t$ with denominator $\prod (1-t^{d_i})$.
(Throughout this note, $\dim$ is vector space dimension, that is,
is $\dim_\C$.)

Now $P(t)$ is known: it was given as a series by Cayley \& Sylvester
(cf. \cite{Sy}) and as a rational function by Springer \cite{Spr}.
For $n=10$ we have
\begin{align*}
P(t) = \phantom{} & 1 + t^2 + 2t^4 + 6t^6 + 12t^8 + 5t^9 + 24t^{10} +
13t^{11} + 52t^{12} + 33t^{13} + 97t^{14} + \\
& 80t^{15} + 177t^{16} + 160t^{17} + 319t^{18} + 301t^{19} + 540t^{20} +
547t^{21} + 887t^{22} + \\
& 926t^{23} + 1429t^{24} + 1512t^{25} + 2219t^{26} + 2402t^{27} +
3367t^{28} + 3681t^{29} + \\
& 5015t^{30} + 5502t^{31} + 7294t^{32} + 8064t^{33} + 10419t^{34} +
11550t^{35} + \\
& 14664t^{36} + 16253t^{37} + 20287t^{38} + 22531t^{39} + 27682t^{40} +
30738t^{41} + \\
& 37319t^{42} + 41378t^{43} + 49671t^{44} + 55060t^{45} + 65390t^{46} +
72391t^{47} + \\
& 85250t^{48} + \cdots
\end{align*}
If we generate invariants of degree $m$, and have found $\dim I_m$
independent ones, then we have found all in degree $m$.
If we know that there is a hsop with degrees
$2$, $4$, $6$, $6$, $8$, $9$, $10$, $14$, then
$$
P(t) = a(t) / (1-t^2)(1-t^4)(1-t^6)^2(1-t^8)(1-t^9)(1-t^{10})(1-t^{14})
$$
where
\begin{align*}
a(t) = \phantom{} & 1 + 2t^6 + 4t^8 + 4t^9 + 7t^{10} + 8t^{11} + 15t^{12} +
15t^{13} + 20t^{14} + 27t^{15} + \\
& 29t^{16} + 35t^{17} + 40t^{18} + 44t^{19} + 47t^{20} + 55t^{21} +
52t^{22} + 57t^{23} + 56t^{24} + \\
& 57t^{25} + 52t^{26} + 55t^{27} + 47t^{28} + 44t^{29} + 40t^{30} +
35t^{31} + 29t^{32} + 27t^{33} + \\
& 20t^{34} + 15t^{35} + 15t^{36} + 8t^{37} + 7t^{38} + 4t^{39} + 4t^{40} +
2t^{42} + t^{48}.
\end{align*}
This means that all basic invariants have degree at most 48,
and we never have to consider subspaces of dimension larger
than 85250, which is doable.

\medskip
So, the procedure is to find basic invariants in some way,
and multiply them together so as to construct for each $m$ the
invariants in $I_m$ that are known already. Compute a basis for
the subspace of $I_m$ spanned by these known invariants,
and if this subspace has the same dimension as $I_m$ itself,
it is all of $I_m$ and we can go to the next $m$.
Since any invariant can be written as a linear combination
of bracket monomials, it seems reasonable to expect that one can
find a spanning set for $I_m$ by just randomly generating some bracket
monomials. This worked fine for the nonic, and for the decimic in degrees
below 21, but in degree 21 where $\dim I_{21} = 547$ and we quickly
generated a subspace of dimension 546, a few dozen attempts to
randomly generate an invariant outside this hyperplane failed.
Therefore, we reverted to the procedure with guaranteed success:
Gordan proved that a basis for the invariants can be found effectively
by computing transvectants, and this indeed yielded the 106th invariant.
(Immediately afterwards the random process also succeeded.)

\medskip
Various reductions simplify the calculations.
First of all, we did the computations modulo a small prime $p$,
e.g. $p=109$ worked. If the images of the invariants under
reduction mod $p$ are independent, then the invariants are independent.
Secondly, if the form is $\sum _{i=0}^{10} \binom{10}{i}a_ix^{10-i}y^i$, we took $a_4=a_7=a_9=0$ and $a_{10}=1$.
Again: if the images of the invariants under this substitution
are independent, then the invariants are independent.
Similar things work for the nonic. But here we have the invariant
$j_2 = a_0a_{10}-10a_1a_9+45a_2a_8-120a_3a_7+210a_4a_6-126a_5^2$ of degree 2.
After the substitutions this becomes $a_0+45a_2a_8-126a_5^2$, and the
substitution $a_0 = -45a_2a_8+126a_5^2$ maps $I_m$ onto $I_m/j_2 I_{m-2}$,
and $\dim I_m/j_2 I_{m-2} = \dim I_m - \dim I_{m-2}$.
Now six variables ($a_1,a_2,a_3,a_5,a_6,a_8$) are left, and the largest dimension
occurring is $\dim I_{48}/j_2I_{46} = 19860$, comparatively small.
(Compared to $\dim I_{48} = 85250$, this saves almost a factor 80
in computation time when an $O(N^3)$ rank algorithm is used.)

\medskip
The computation was done, and the result is: for $m \le 48$ the
values of $d_m$ are as listed in Proposition \ref{prop1}.
Consequently, if there is a system of parameters with degrees
2, 4, 6, 6, 8, 9, 10, 14, so that no basic invariant has degree
larger than 48, then Proposition \ref{prop1} follows.

\section{A system of parameters for $\oo (V_{10})^{\SLtwo}$}
Let $V_n$ be the space of forms of degree $n$ (in the variables $x,y$).
A {\it covariant of order m and degree d} of $V_n$ is an
$\SLtwo$-equivariant homogeneous polynomial map
$\phi: V_n \rightarrow V_m$ of degree $d$.
The invariants of $V_n$ are the covariants of order 0.
The identity map is a covariant of order $n$ and degree 1.
Customarily, one indicates such a covariant $\phi$ by giving its
image of a generic element $f \in V_n$.
(In particular, the identity map is noted $f$.)
Let $V_{m,d}$ be the space of covariants of order $m$ and degree $d$.

Consider $f \in V_{10}$, 
\[f=a_0x^{10}+10a_1x^9y+\ldots +10 a_9xy^9+a_{10}y^{10}, \]
and the following covariants
\begin{alignat*}{4}
k&=(f,f)_8 &\in\,& V_{4,2},\quad\quad&
m&=(f,k)_4 &\in\,& V_{6,3},\quad\quad\quad\quad\\
q&=(f,f)_6 &\in\,& V_{8,2},&
r&=(f,q)_8 &\in\,& V_{2,3},\\
k_q&=(q,q)_6 &\in\,& V_{4,4},&
k_m&=(m,m)_4\, &\in\,& V_{4,6},\\
m_q&=(q,k_q)_4~ &\in\,& V_{4,6},
\end{alignat*}
and invariants (the suffix indicates the degree)
\begin{alignat*}{4}
j&_2 &\,=\,& (f,f)_{10},&
A&_6 &\,=\,& (m,m)_6,\\
j&_4 &\,=\,& (k,k)_4,&
C&_6 &\,=\,& (r,r)_2,\\
j&_8 &\,=\,& (k,k_m)_4,&
j&_{14} &\,=\,& ((k_q,k_q)_2,m_q)_4,\\
j&_9 &\,=\,& ((m,k)_1,k^2)_8,\quad\quad&
A&_{14} &\,=\,& ({(k,k)_2}^2,(m,m)_2)_8,\\
j&_{10} &\,=\,& ((m,m)_2,k^2)_8.
\end{alignat*}

\smallskip
\begin{theorem}\label{explhsop}
The eight invariants $j_2$, $j_4$, $A_6$, $C_6$, $j_8$, $j_9$,
$j_{10}$, $j_{14}+A_{14}$ form a homogeneous system of parameters for
the ring $\oo(V_{10})^{\SLtwo}$ of invariants of the binary decimic.
\end{theorem}

This is proved by invoking Hilbert's characterization of
homogeneous systems of parameters as sets that define the nullcone.

\subsection{The nullcone}

The {\em nullcone} of $V_n$, denoted $\nn(V_n)$, is the set of binary forms
of degree $n$ on which all invariants vanish. It turns out (\cite{Hi2})
that this is precisely the set of binary forms of degree $n$
with a root of multiplicity $>\frac{n}{2}$. The elements of $\nn(V_n)$
are called {\em nullforms}.
The nullcone $\nn(V_n\oplus V_m)$ is
the set of pairs $(g,h)\in V_n\oplus V_m$ such that $g$ and $h$ have a
common root of multiplicity $>\frac{n}{2}$ in $g$ and of multiplicity
$>\frac{m}{2}$ in $h$. (In this note, this result can be taken as the
definition of the symbol $\nn(V_n\oplus V_m)$.)

We have the following result, due to Hilbert \cite{Hi2},
formulated for the particular case of binary forms:
\begin{proposition} \label{hilbert} 
For $n \ge 3$, consider $i_1,\ldots ,i_{n-2}\in \oo(V_n)^{\SLtwo}$
homogeneous invariants of $V_n$.
The following two conditions are equivalent:
\begin{itemize}
\item[(i)] $\nn(V_n)=\vv(i_1,\ldots ,i_{n-2})$,
\item[(ii)] $\{i_1,\ldots ,i_{n-2}\}$ is a homogeneous system of parameters
of $\oo(V_n)^{\SLtwo}$. 
\end{itemize}
\end{proposition}
(Here $\vv(J)$ stands for the vanishing locus of $J$.)

\medskip\noindent
We prove the above theorem by first finding a defining set
for the nullcone that is still too large, and then showing
that some elements are superfluous.

\medskip
\medskip
We need information on the invariants of $V_n$ for $n = 2,\,4,\,6,\,8$:
\begin{lemma}\label{hsops}
The following are systems of parameters of
$\oo(V_n)^{\SLtwo}$ for $n = 2,\,4,\,6,\,8$.
\begin{itemize}
\item[(i)] If $n=2$: $(f,f)_2$ of degree $2$.
\item[(ii)] If $n=4$: $(f,f)_4$ and $((f,f)_2,f)_4$ of degrees $2$ and $3$.
\item[(iii)] If $n=6$: $(f,f)_6$, $(k,k)_4$, $((k,k)_2,k)_4$,
and $(m^2,(k,k)_2)_4$ of degrees $2$, $4$, $6$ and $10$,
where $k = (f,f)_4$ and $m = (f,k)_4$.
\item[(iv)] If $n=8$: $(f,f)_8$, $((f,f)_4,f)_8$, $(k,k)_4$, $(m,k)_4$,
$((k,k)_2,k)_4$, $((k,k)_2,m)_4$ of degrees $2$, $3$, $4$, $5$, $6$ and $7$,
where $k = (f,f)_6$ and $m = (f,k)_4$.
\end{itemize}
\end{lemma}
\begin{proof}
This is classical for $n=2$, $4$, $6$, see, e.g., \cite{Clebsch,GrYo,Schur},
and due to von Gall \cite{vG8} and Shioda \cite{Sh} for $n=8$.
\end{proof}

\begin{lemma} \label{jerzy} {\rm (Weyman \cite{We1})}
Let $f\in V_d$. If $d > 4k-4$ and all $(f,f)_{2k}$, $(f,f)_{2k+2}$, ...
vanish, then $f$ has a root of multiplicity $d-k+1$. If $d = 4k-4$ and
$((f,f)_{2k-2},f)_d$, $(f,f)_{2k}$, $(f,f)_{2k+2}$, ... vanish,
then $f$ has a root of multiplicity $d-k+1$. \qed
\end{lemma}

\begin{lemma} \label{covariants} Let $f\in V_{10}$ and $j_2=(f,f)_{10}$,
$k=(f,f)_8\in V_4$, $m=(f,k)_4 \in V_6$, $q=(f,f)_6\in V_8$. We have:
\begin{enumerate}
\item[(i)] If $j_2=0$, $k\neq 0$ and $(k,m)\in \nn_{V_4\oplus V_6}$,
then $f$ has a root of multiplicity 6.
\item[(ii)] If $j_2=0$, $k=0$ and $0\neq q \in \nn_{V_8}$,
then $f$ has a root of multiplicity 7.
\item[(iii)]
If $j_2=0$, $k=0$ and $q=0$, then $f$ has a root of multiplicity 8.  
\end{enumerate}
\end{lemma}
\begin{proof}
The covariants $k$, $q$ and the invariant $j_2$ are:
\begin{align*}
j_2&=-252 a_5\supr{2} + 420 a_4 a_6 - 240 a_3 a_7 + 90 a_2 a_8 - 20 a_1 a_9 + 2 a_0 a_{10},\\
k&=(70 a_6\supr{2} - 112 a_5 a_7 + 56 a_4 a_8 - 16 a_3 a_9 + 2 a_2 a_{10})y^4+\\
&(56 a_5 a_6 - 112 a_4 a_7 + 80 a_3 a_8 - 28 a_2 a_9 + 4 a_1 a_{10}) xy^3+\\
&(168 a_5\supr{2} - 252 a_4 a_6 + 96 a_3 a_7 - 6 a_2 a_8 -8 a_1 a_9 + 2 a_0 a_{10}) x^2y^2+\\
&(56 a_4 a_5 - 112 a_3 a_6 + 80 a_2 a_7 - 28 a_1 a_8 + 4 a_0 a_9) x^3 y+\\
&(70 a_4\supr{2} - 112 a_3 a_5 + 56 a_2 a_6 - 16 a_1 a_7 + 2 a_0 a_8)x^4,\\
q&=(-20 a_7\supr{2} + 30 a_6 a_8 - 12 a_5 a_9 +  2 a_4 a_{10})y^8+\\
&(-40 a_6 a_7 + 72 a_5 a_8 - 40 a_4 a_9 + 8 a_3 a_{10})y^7 x+\\
&( -140 a_6\supr{2} + 168 a_5 a_7 - 40 a_3 a_9 + 12 a_2 a_{10})y^6x\supr{2}+\\
&( -168 a_5 a_6 + 280 a_4 a_7 - 120 a_3 a_8 + 
 8 a_1 a_{10}) y^5x^3+\\
 &(-252 a_5\supr{2} + 280 a_4 a_6 + 40 a_3 a_7 - 
 90 a_2 a_8 + 20 a_1 a_9 + 2 a_0 a_{10}) y^4x^4+\\
 &(-168 a_4 a_5 + 
 280 a_3 a_6 - 120 a_2 a_7 + 8 a_0 a_9) y^3x^5+\\
 &(-140 a_4\supr{2} + 
 168 a_3 a_5 - 40 a_1 a_7 + 12 a_0 a_8) y^2x^6+\\
 &(-40 a_3 a_4 + 
 72 a_2 a_5 - 40 a_1 a_6 + 8 a_0 a_7) y x^7+\\
 &(-20 a_3\supr{2} + 
 30 a_2 a_4 - 12 a_1 a_5 + 2 a_0 a_6)x^8.
\end{align*}

\noindent
(i). If $(k,m) \in \nn_{V_4\oplus V_6}$ then $k$ and $m$ have a common root,
of multiplicity 3 in $k$ and of multiplicity 4 in $m$.
Without loss of generality we consider the cases $k=x^4$, $x^4 \divides m$
and $k=x^3y$, $x^4 \divides m$. 

\medskip\noindent
Case 1: $k=x^4$. Then $m$ becomes:
\[m=(f,x^4)_4=a_4 x^6 + 6 a_5 x^5 y + 15 a_6 x^4 y^2 + 20 a_7 x^3 y^3 + 
 15 a_8 x^2 y^4 + 6 a_9 x y^5 + a_{10} y^6.\]
From $x^4 \divides m$ it follows $a_7=\ldots =a_{10}=0$. We replace this in $k$ and because we supposed $k=x^4$ we obtain also $a_6=a_5=0$. But then $x^6 \divides f$, hence $f$ will have a root of multiplicity 6. 

\medskip\noindent 
Case 2: $k=x^3 y$. Then $m$ becomes:
\[
m\!=\!(f,x^3 y)_4\!=\!
-a_3 x^6 - 6 a_4 x^5 y - 15 a_5 x^4 y^2 - 20 a_6 x^3 y^3 - 
 15 a_7 x^2 y^4 - 6 a_8 x y^5 - a_9 y^6.
\]
From $x^4 \divides m$ it follows $a_6=\ldots =a_9=0$. We replace this in $k$ and $j_2$ and as we supposed $k=x^3y$ we obtain 
\begin{align*}
168 a_5\supr{2} + 2 a_0 a_{10}=0,\\
-252 a_5\supr{2} + 2 a_0 a_{10}=0,
\end{align*}
which implies $a_5=0$. But then the coefficient of $x^3$ in $k$ becomes 0. Contradiction with our assumption. 

\medskip\noindent
(ii). Without loss of generality we suppose $x^5 \divides q$.
We denote by $J$ the ideal generated by $j_2$, the coefficients of $k$
and the coefficients of $x^4 y^4,x^3y^5,\ldots, y^8$ in $q$.
Denote also by $p_1$, $p_2$ and $p_3$ the coefficients of
$x^7y$, $x^6y^2$ and $x^5y^3$, respectively, in $q$. We have 
\[p_1\supr{4},p_2\supr{3},p_3\supr{2} \in J,\]   
which implies that $x^8 \divides q$. 

Consider now the ideal $J$ generated by $j_2$, the coefficients of $k$
and the coefficients of $x^7y,x^6y^2,\ldots, y^8$ in $q$.
Denote by $p_0$ the coefficient of $x^8$ in $q$.
We have $a_i p_0\in J$ for $i=10,9,8,7,6,5,4$.
Because $q \neq 0$ we find $a_{10}=\ldots=a_4=0$.
This means that $x^7 \divides f$, so $f$ will have a root of multiplicity 7.

\medskip\noindent
(iii). This follows from Lemma \ref{jerzy}. 
\end{proof}

\begin{lemma} \label{kmv10}
Let $k\in V_4$ and $m\in V_6$, $k\neq 0$, $m\neq 0$, both of them nullforms.
If the transvectants $((m,m)_4,k)_4$, $((m,m)_2,k^2)_8$, $(m^2,k^3)_{12}$,
$((m,k)_1,k^2)_8$, and $((k,k)_2\supr{2},(m,m)_2)_8$ vanish, then
$(k,m)\in \nn_{V_4\oplus V_6}$.
\end{lemma}
\begin{proof}
Suppose $(k,m)\notin \nn_{V_4\oplus V_6}$.
Without loss of generality we suppose
\begin{align*}
k=&x^3(a_1 x+a_2y),\\
m=&y^4(b_1x^2+b_2xy+b_3y^2).
\end{align*}
We have
\[0=((m,m)_4,k)_4 \sim a_1b_1\supr{2}\] 

\medskip\noindent
Case 1: $a_1=0$. Then 
\begin{alignat*}{2}
0&=((m,m)_2,k^2)_8 &\,\sim\,& a_2\supr{2} b_1\supr{2},\\ 
0&=((m,k)_1,k^2)_8 &\,\sim\,& a_2\supr{3}b_3,\\ 
0&=((k,k)_2\supr{2},(m,m)_2)_8 &\,\sim\,& a_2\supr{4}(5b_2\supr{2}-12 b_1b_3) 
\end{alignat*}
Because $k\neq 0$ we have $a_2 \neq 0$, but then it follows that
$b_1=b_3=b_2=0$. Contradiction with $m\neq 0$.

\medskip\noindent
Case 2: $a_1 \ne 0$, $b_1=0$. Then 
\begin{alignat*}{2}
0&=((m,m)_2,k^2)_8 &\,\sim\,& a_1\supr{2} b_2\supr{2},\\ 
0&=((m,k)_1,k^2)_8 &\,\sim\,& a_2\supr{3}b_3,\\ 
0&=((k,k)_2\supr{2},(m,m)_2)_8 &\,\sim\,& a_2\supr{4}b_2\supr{2},\\ 
0&=(m^2,k^3)_{12} &\,\sim\,& a_1(a_2\supr{2} b_2\supr{2} - 11 a_1 a_2 b_2 b_3
+ 22 a_1\supr{2} b_3\supr{2}) 
\end{alignat*}
If $a_2\neq 0$ then $b_2=b_3=0$.
And if $a_2=0$ then $a_1\supr{2}b_2\supr{2}=a_1\supr{3} b_3\supr{2}=0$,
and again $b_2=b_3=0$. Contradiction with $m\neq 0$.
\end{proof}

%
%
%

After this preparation we can write down a defining set for the nullcone.
Define $k$, $m$, $q$, $j_2$, $j_4$, $A_6$, $j_8$, $j_9$, $j_{10}$, $j_{14}$,
$A_{14}$ as above (before Theorem \ref{explhsop}), and moreover
\begin{alignat*}{4}
j&_6 &\,=\,& ((k,k)_2,k)_4,\quad\quad&
A&_{12} &\,=\,& (m^2,k^3)_{12},\\
B&_6 &\,=\,& ((q,q)_4,q)_8.
\end{alignat*}

\begin{proposition} \label{nullconev10}
With notations as above, the nullcone $\nn _{V_{10}}$ is defined by
\[
\nn _{V_{10}}=\vv (j_2,j_4,j_6,A_6,B_6,j_8,j_9,j_{10},A_{12},j_{14},A_{14}).
\]
\end{proposition}
\begin{proof}
Since $k \in V_4$ we can apply Lemma \ref{hsops}(ii) and conclude that
if $j_4=j_6=0$ then $k$ is a nullform.
Without loss of generality we consider three cases:
$k=0$, $k=x^4$ and $k=x^3y$. 

\medskip\noindent
Case 1: $k=0$. Denote by $I=(j_2,k)$ the ideal generated by $j_2$ and
the coefficients of $k$. Define
\begin{alignat*}{4}
A&_4 &\,=\,& (q,q)_8,&
A&_{10} &\,=\,& (m_q,k_q)_4,\\
A&_8 &\,=\,& (k_q,k_q)_4,\quad\quad&
B&_{12} &\,=\,& ((k_q,k_q)_2,k_q)_4,
\end{alignat*}
Since $q \in V_8$, in order to show that $q$ is a nullform
it suffices by Lemma \ref{hsops}(iv) to show that each of
$A_4$, $B_6$, $A_8$, $A_{10}$, $B_{12}$ and $j_{14}$ vanishes.

Easy Gr\"obner basis computations show that
$A_4,A_8,A_{10} \in I$ and $B_{12} \in (I,B_6)$.
It follows that if $k=0$ and $j_2=B_6=j_{14}=0$ then $q$ is a nullform.
Now Lemma \ref{covariants} implies that $f$ is a nullform. 

\medskip\noindent
Case 2: $k=x^4$. Then we have: 
\begin{align*}
A_{12} &\sim a_{10}\supr{2},\\
j_{10} &\sim - a_9\supr{2} + a_8 a_{10},\\
j_8 &\sim 3 a_8\supr{2} - 4 a_7 a_9 + a_6 a_{10},\\
A_6 &\sim -10 a_7\supr{2} + 15 a_6 a_8 - 6 a_5 a_9 + a_4 a_{10}.
\end{align*}  
If $A_{12}=j_{10}=j_8=A_6=0$ then it follows that $a_{10}=\ldots =a_7=0$.
If we substitute this in $k$ we obtain
\begin{align*}
k=\,&70 a_6\supr{2} y^4+56 a_5 a_6  xy^3+
(168 a_5\supr{2} - 252 a_4 a_6) x^2y^2+\\
&(56 a_4 a_5 - 112 a_3 a_6) x^3 y+
(70 a_4\supr{2} - 112 a_3 a_5 + 56 a_2 a_6)x^4,
\end{align*}
and as we supposed $k=x^4$ we get also $a_6=a_5=0$,
which implies that $f$ is a nullform. 

\medskip\noindent 
Case 3: $k=x^3y$. Then we have:
\begin{align*}
j_9 &\sim a_9,\\ 
A_{14} &\sim a_7a_9-a_8\supr{2},\\ 
j_{10} &\sim -5 a_7\supr{2} + 2 a_6 a_8 + 3 a_5 a_9,\\ 
A_6 &\sim -10 a_6\supr{2} + 15 a_5 a_7 - 6 a_4 a_8 + a_3 a_9. 
\end{align*}
If $j_9=A_{14}=j_{10}=A_6=0$ then $a_9=\ldots =a_6=0$.
We substitute this in $k$ and $j_2$:
\begin{align*}
k=\,& 2 a_2 a_{10}y^4+4 a_1 a_{10} xy^3+
(168 a_5\supr{2}  + 2 a_0 a_{10}) x^2y^2+\\
&56 a_4 a_5  x^3 y+
(70 a_4\supr{2} - 112 a_3 a_5 )x^4,\\
j_2=\,&-252 a_5\supr{2} + 2 a_0 a_{10}
\end{align*}
From $168 a_5\supr{2}  + 2 a_0 a_{10}=-252 a_5\supr{2} + 2 a_0 a_{10}=0$
we find $a_5=0$, which contradicts $k=x^3y$. 
\end{proof}

So far, we defined the nullcone using 11 invariants, but we need
a definition using 8 invariants. As a first step, replace the two
invariants of degree 14 by a single one.

Now for $f=x^2 y (2a_1 x^7 + 9a_8 y^7)$
all invariants from Proposition \ref{nullconev10} vanish, except $A_{14}$.
And for $f=y^3 (120a_3 x^7 + a_{10} y^7)$
all invariants from Proposition \ref{nullconev10} vanish, except $j_{14}$.
That means that the single invariant of degree 14 cannot be either
$j_{14}$ or $A_{14}$. However, as it turns out we can use $j_{14}+A_{14}$.

\subsection{Finding the system of parameters}
Proposition \ref{nullconev10} gives an explicit set of invariants
(and in particular an explicit set of degrees of invariants)
that define the nullcone. Having that, only a finite amount of work is left.

The final part of the construction of the system of parameters
was done by computer. All computations were carried out in the
ring $R$ generated by the 106 invariants found in Section \ref{finding}.
Or, more precisely, in the quotient $Q = R/j_2R$, reduced mod $p$,
where this time $p = 197$ (the different $p$ has no significance),
and again $a_4$, $a_7$ and $a_9$ were taken to be zero.
It was checked that the graded parts of the resulting ring have
the expected dimension (for degree up to 54), so that no collapse
occurred as a consequence of the reduction mod $p$ or the substitution
of variables.

The ideal generated in this ring by all invariants of degrees
4, 6, 8, 9, 10, 14 has full dimension 542 for its graded part
of degree 24. We know that $\dim I_{24} = 1429$ and $\dim I_{22} = 887$
and multiplication by $j_2$ is an injection, so $\dim I_{24}/j_2I_{22} = 542$.
It follows that the ideal generated by these invariants, together with $j_2$,
contains all of $I_{24}$, so that no invariants of degree 12 are needed
to define the nullcone (since their squares are in $I_{24}$, and they
themselves are in the radical).

With only $j_{14}+A_{14}$ instead of all invariants of degree 14
in the set of generators of the ideal, one finds full dimension 1148
for the graded part of degree 28, so this single invariant of degree 14
suffices.

With only $j_{10}$ instead of all invariants of degree 10, one finds
full dimension 221 in degree 20, so this single invariant of degree 10
suffices.

With only $j_9$ instead of all invariants of degree 9, one finds
full dimension 890 in degree 27, so this single invariant of degree 9
suffices.

With only $j_8$ instead of all invariants of degree 8, one finds
full dimension 2279 in degree 32, so this single invariant of degree 8
suffices.

That only leaves the invariants of degree 6. After some work it turned out
that with only $A_6$ and $C_6$ one finds full dimension 37892 in degree 54,
so these suffice, and we have constructed the homogeneous system of
parameters promised in Theorem \ref{explhsop}.

\medskip
Note that one knows what to expect if all is well:
the coefficients of the polynomial $a(t)$ from Section \ref{finding}
give for each degree the codimension of the set of invariants
in the ideal generated by the hsop in the space of all invariants
of that degree. Since 54 is the smallest multiple of 6 where $a(t)$
has zero coefficient, that explains why the computation had to extend
to there.

\medskip\noindent
Addresses of authors:

\smallskip
\begin{minipage}{2in}
Andries E. Brouwer \\
Dept. of Math. \\
Techn. Univ. Eindhoven \\
P. O. Box 513 \\
5600MB Eindhoven \\
Netherlands \\
{\tt aeb@cwi.nl}
\end{minipage}
\begin{minipage}{2in}
Mihaela Popoviciu \\
Mathematisches Institut \\
Universit\"at Basel \\
Rheinsprung 21 \\
CH-4051 Basel \\
Switzerland \\
{\tt mihaela.popoviciu@unibas.ch}
\end{minipage}

\end{document}